\documentclass{article}

\usepackage{amsmath}
\usepackage{amssymb}

\newtheorem{thm}{Theorem}
\newtheorem{conj}{Conjecture}
\newtheorem{lem}{Lemma}

\title{A note on powerful numbers in short intervals}
\author{Tsz Ho Chan}
\date{}

\begin{document}
\maketitle

\begin{abstract}
In this note, we are interested in obtaining uniform upper bounds for the number of powerful numbers in short intervals $(x, x + y]$. We obtain unconditional upper bounds $O(\frac{y}{\log y})$ and $O(y^{11/12})$ for all powerful numbers and $y^{1/2}$-smooth powerful numbers respectively. Conditional on the $abc$-conjecture, we prove the bound $O(\frac{y}{\log^{1+\epsilon} y})$ for squarefull numbers and the bound $O(y^{(2 + \epsilon)/k})$ for $k$-full numbers when $k \ge 3$. They are related to Roth's theorem on arithmetic progressions and the conjecture on non-existence of three consecutive squarefull numbers.
\end{abstract}

\section{Introduction and Main Result}

A number $n$ is {\it squarefull} if its prime factorization $n = p_1^{a_1} p_2^{a_2} \cdots p_r^{a_r}$ satisfies $a_i \ge 2$ for all $1 \le i \le r$. Similarly, a number $n$ is $k$-{\it full} if the exponents $a_i \ge k$ for all $1 \le i \le r$. For example, $72 = 2^3 \cdot 3^2$ is squarefull and $243 = 3^5$ is $5$-full. Let $Q_k(x)$ denote the number of $k$-full numbers which are less than or equal to $x$. It is known that
\begin{equation} \label{asymptotic}
Q_k(x) = \prod_{p} \Bigl(1 + \sum_{m = k+1}^{2k - 1} \frac{1}{p^{m/h}} \Bigr) x^{1/k} + O(x^{1/(k+1)})
\end{equation}
where the product is over all primes (see \cite{ES} or \cite{BG} for example.) There are also interests in counting the number of $k$-full numbers in short intervals $(x, x+y]$ with $y = o(x)$. For moderate size $y$, there are some asymptotic results: For example, Trifonov \cite{T} and Liu \cite{L} obtained,
\begin{equation} \label{shortsquare}
Q_2 \bigl(x + x^{1/2 + \theta} \bigr) - Q_2(x) \sim \frac{\zeta(3/2)}{2 \zeta(3)} x^\theta \; \; \text{ for } \; \; 19 / 154 = 0.12337... < \theta < 1/2,
\end{equation}
and
\begin{equation} \label{shortcube}
Q_3 \bigl(x + x^{2/3 + \theta} \bigr) - Q_3(x) \sim \frac{\zeta(4/3)}{3 \zeta(4)} x^\theta \; \; \text{ for } \; \; 5/42 = 0.11904... < \theta < 1/3
\end{equation}
respectively.

\bigskip

What happens when $y$ is very small such as $y \ll x^{1/2}$ or even $y \ll \log x$? For such short intervals, one can only expect suitable upper bounds rather than asymptotic formulas. Thus, in this note, we are interested in finding uniform upper bounds for $Q_k(x + y) - Q_k(x)$ with $1 \le y \le x$ that are independent of $x$. By comparing $k$-full numbers with perfect $k$-th powers, we suspect the following to be true.
\begin{conj} \label{mainconj}
Given an integer $k \ge 2$ and a real number $x \ge 1$. There exists some constant $C_k \ge 1$ such that
\[
Q_k(x + y) - Q_k(x) \le C_k y^{1/k}
\]
uniformly over $1 \le y \le x$.
\end{conj}
We are far from proving this at the moment. The current best upper bound
\begin{equation} \label{upper1}
Q_k(x + y) - Q_k(x) \ll \frac{y \log \log (y + 2)}{\log (y + 2)}
\end{equation}
was obtained by De Koninck, Luca and Shparlinski \cite{DLS}. In this note, we improve \eqref{upper1} slightly.
\begin{thm} \label{thm1}
Given an integer $k \ge 2$ and a real number $x \ge 1$. We have
\begin{equation} \label{upper2}
Q_k(x + y) - Q_k(x) \ll \frac{y}{\log (y + 1)}
\end{equation}
uniformly over $1 \le y \le x$.
\end{thm}
In fact, we shall prove the following more general result concerning squarefull numbers in arithmetic progressions over short intervals which gives Theorem \ref{thm1} immediately as $k$-full numbers are included in squarefull numbers.
\begin{thm} \label{thm2}
Given real numbers $x \ge 1$ and $0 < \alpha < 1$. Suppose $q > 0$ and $r$ are two integers with $\text{gcd}(r, q) = 1$. We have
\begin{equation*}
\mathop{\mathop{\sum_{x < n \le x + y}}_{n \text{ squarefull}}}_{n \equiv r \;(\bmod{q})} 1 \ll_\alpha \frac{y}{\phi(q) \log (y + 1)}
\end{equation*}
uniformly over $1 \le y \le x$ and $1 \le q \le y^{1 - \alpha}$.
\end{thm}
Using similar technique, one can obtain some power savings over \eqref{upper2} for {\it smooth} $k$-full numbers in short intervals.
\begin{thm} \label{thm1.5}
Given an integer $k \ge 2$ and a real number $x \ge 1$. We have
\begin{equation} \label{upper1.5}
\mathop{\mathop{\sum_{x < n \le x + y}}_{n \text{ $k$-full}}}_{p^{+}(n) \le y^{1/2}} 1 \le \mathop{\mathop{\sum_{x < n \le x + y}}_{n \text{ squarefull}}}_{p^{+}(n) \le y^{1/2}} 1 \ll y^{11/12}.
\end{equation}
uniformly over $1 \le y \le x$. Here $p^{+}(n)$ stands for the largest prime factor of $n$. Remark: One may increase the exponent $1/2$ up to $1$ and obtain similar power saving upper bound.
\end{thm}

The bound \eqref{upper1.5} lends evidence towards Conjecture \ref{mainconj} and shows that the difficulty lies with {\it non-smooth} $k$-full numbers. Another piece of evidence comes from the famous $abc$-conjecture. It was proved in \cite{DLS} that, given any $\epsilon > 0$, the interval
\begin{equation} \label{interval}
(x, x + x^{1 - (2 + \epsilon) / k}]
\end{equation}
contains at most one $k$-full number for sufficiently large $x$ under the $abc$-conjecture. From this, one has
\begin{thm} \label{thm2.5}
Given an integer $k \ge 3$ and real numbers $\epsilon > 0$ and $x \ge 1$. We have
\begin{equation} \label{upper3}
Q_k(x + y) - Q_k(x) \ll_\epsilon y^{(2 + \epsilon) / k}
\end{equation}
uniformly over $1 \le y \le x$.
\end{thm}

Unfortunately, the proof in \cite{DLS} concerning \eqref{interval} was a little inaccurate as the $a$, $b$, $c$ in its application of the $abc$-conjecture might not be relatively prime. So, we shall give a complete proof of Theorem \ref{thm2.5} later. Furthermore, \eqref{interval} gives us nothing when $k = 2$. To remedy this, we shall prove the following conditional result which improves on \eqref{upper2} slightly by a small power of logarithm.
\begin{thm} \label{thm3}
The $abc$-conjecture implies that, for some absolute constant $c > 0$,
\[
Q_2(x+y) - Q_2(x) \ll \frac{y}{\log^{1+c} (y+1)}
\]
uniformly over $1 \le y \le x$.
\end{thm}
Its proof relies on the following recent breakthrough result of Bloom and Sisask \cite{BS} on density of integer sequence without three term arithmetic progressions.
\begin{thm}[Bloom-Sisask] \label{BSthm}
Let $N \ge 2$ and $A \subset \{1, 2, ..., N \}$ be a set with no non-trivial three-term arithmetic progressions, i.e. solutions to $x + y = 2z$ with $x \neq y$. Then
\[
|A| \ll \frac{N}{(\log N)^{1+c}},
\]
where $c > 0$ is an absolute constant.
\end{thm}

\bigskip

This paper is organized as follows. First, we will prove Theorems \ref{thm2} and \ref{thm1.5} using Brun-Titchmarsh inequality and ideas from Shiu's generalization \cite{S}. Then we will prove Theorem \ref{thm2.5} using the $abc$-conjecture. Finally, we will prove Theorem \ref{thm3} by establishing non-existence of three term arithmetic progressions for squarefull numbers in short intervals.

\bigskip

{\bf Notation.} We use $| A |$ to denote the number of elements in a finite set $A$ and $\lfloor x \rfloor$ to denote the greatest integer less than or equal to $x$. We let $p_{-}(n)$ and $p^{+}(n)$ be the smallest and the largest prime factor of $n$ respectively. The symbols $f(x) = O(g(x))$, $f(x) \ll g(x)$ and $g(x) \gg f(x)$ are equivalent to $|f(x)| \leq C g(x)$ for some constant $C > 0$. Also, $f(x) = O_{\lambda} (g(x))$, $f(x) \ll_{\lambda} g(x)$ or $g(x) \gg_{\lambda} f(x)$ mean that the implicit constant may depend on $\lambda$. Furthermore, $f(x) = o(g(x))$ means $\lim_{x \rightarrow \infty} f(x)/g(x) = 0$ and $f(x) \sim g(x)$ means $\lim_{x \rightarrow \infty} f(x)/g(x) = 1$. Finally, the summation symbol $\mathop{\sum\nolimits'}$ signifies that a sum is over squarefull numbers only.


\section{Some Preparations}

\begin{lem} \label{lem1}
For any $X \ge 1$,
\[
\mathop{\sum\nolimits'}_{X < n \le X^2} \frac{1}{n} \ll X^{-1/2}.
\]
\end{lem}

Proof: From \eqref{asymptotic}, we have $Q_2(X) \ll X^{1/2}$. By partial summation, the above sum is
\[
= \int_{X}^{X^2} \frac{1}{u} d Q(u) = \frac{Q(X^2)}{X^2} - \frac{Q(X)}{X} + \int_{X}^{X^2} \frac{Q(u)}{u^2} du \ll \frac{1}{X^{1/2}} + \int_{X}^{X^2} \frac{1}{u^{3/2}} du \ll \frac{1}{X^{1/2}}.
\]
\begin{lem}[Brun-Titchmarsh inequality] \label{lem2}
Let $q \ge 1$ and $r$ be integers satisfying $\text{gcd}(r, q) = 1$. Suppose $q < y \le x$ and $z \ge 2$. Then
\[
\mathop{\mathop{\sum_{x < n \le x + y}}_{n \equiv r \; (\bmod{q})}}_{p_{-}(n) > z} 1 \ll \frac{y}{\phi(q) \log z} + z^2.
\]
Remark: The above bound is still true when $y \le q$ or $y < 1$ since there is at most one term in the sum.
\end{lem}

Proof: This follows from Selberg upper bound sieve method (see \cite{HR}, page 104 for example).

\bigskip

Finally, let us recall the $abc$-conjecture. For any nonzero integer $m$, let
\[
\kappa(m) := \prod_{p | m} p
\]
be the kernel of $m$.
\begin{conj}[$abc$-conjecture] \label{abc}
For any $\epsilon > 0$, there exists a constant $C_\epsilon > 0$ such that, for any integers $a, b, c$ with $a + b = c$ and $\text{gcd}(a,b) = 1$, the bound
\[
\max\{ |a|, |b|, |c| \} \le C_\epsilon \kappa(a b c)^{1 + \epsilon}
\]
holds.
\end{conj}


\section{Proof of Theorem \ref{thm2}}

Our proof is inspired by Shiu \cite{S} on Brun-Titchmarsh theorem for multiplicative functions. We may assume that $y \ge 2^{2/\alpha}$ for the theorem is clearly true when $1 \le y < 2^{2/\alpha}$ by choosing a large enough implicit constant. Recall that $1 \le q \le y^{1 - \alpha}$ for some $\alpha > 0$. Let $z = y^{\alpha / 2} \ge 2$. Any squarefull number $n$ in $[x, x + y]$ can be factored as
\[
n = \underbrace{p_1^{a_1} \cdots p_j^{a_j}}_{b_n} \underbrace{p_{j+1}^{a_{j+1}} \cdots p_s^{a_s}}_{d_n} \; \; \text{ with } \; \; p_1 < p_2 < \cdots < p_s
\]
where $j$ is the greatest index such that $ p_1^{a_1} \cdots p_j^{a_j} \le z$. Hence, $b_n \le z < b_n p_{j+1}^{a_{j+1}}$. Note that $j$ may be $0$ (the product is an empty product) if $p_1^{a_1} > z$. In this case, $b_n = 1$ and $d_n = n$. Also, since $n \equiv r \pmod{q}$ with $\text{gcd}(r,q) = 1$, we must have $\text{gcd}(b_n, q) = 1 = \text{gcd}(d_n, q)$.

\bigskip

Case 1: $b_n > z^{1/2}$. As $q \le y^{1 - \alpha}$ and $z = y^{\alpha/2}$, the number of such squarefull numbers is bounded by
\begin{equation} \label{c1}
\mathop{\mathop{\sum\nolimits'}_{z^{1/2} < b \le z}}_{\text{gcd}(b, q) = 1} \mathop{\mathop{\sum_{x < n \le x + y}}_{b | n}}_{n \equiv r \;(\bmod{q})} 1 \le  \mathop{\sum\nolimits'}_{z^{1/2} < b \le z} \Bigl( \frac{y/b}{q} + 1 \Bigr) \ll \frac{y}{q z^{1/4}} + z^{1/2} \ll_\alpha \frac{y}{\phi(q) \log y}
\end{equation}
by \eqref{asymptotic} and Lemma \ref{lem1}.

\bigskip

Case 2: $b_n \le z^{1/2}$ and $p_{-}(d_n) \le z^{1/2}$. Then $p_{j+1} \le z^{1/2}$ and $p_{j+1}^{a_{j+1}} > z^{1/2}$ which implies $p_{j+1}^{-a_{j+1}} \le \min(z^{-1/2}, p_{j+1}^{-2})$ as $a_{j+1} \ge 2$. Hence, the sum
\[
\sum_{p_{j+1} \le z^{1/2}} \frac{1}{p_{j+1}^{a_{j+1}}} \le \sum_{p_{j+1} \le z^{1/4}} z^{-1/2} + \sum_{z^{1/4} < p_{j+1} \le z^{1/2}} \frac{1}{p_{j+1}^2} \ll \frac{1}{z^{1/4}}.
\]
Therefore, by replacing $p_{j+1}^{a_{j+1}}$ with a generic $p^a$, the number of squarefull numbers in this case is bounded by 
\begin{equation} \label{c2}
\mathop{\sum_{p \le z^{1/2}}}_{\text{gcd}(p,q) = 1} \mathop{\mathop{\sum_{x < n \le x + y}}_{p^a | n}}_{n \equiv r \;(\bmod{q})} 1 \le \sum_{p \le z^{1/2}} \Bigl(\frac{y / p^a}{q} + 1 \Bigr) \ll \frac{y}{q z^{1/4}} + z^{1/2} \ll_\alpha \frac{y}{\phi(q) \log y}
\end{equation}
since $q \le y^{1 - \alpha}$ and $z = y^{\alpha/2}$.

\bigskip

Case 3: $b_n \le z^{1/2}$ and $p_{-}(d_n) > z^{1/2}$. As $q \le y^{1 - \alpha}$ and $z = y^{\alpha/2}$, the number of such squarefull numbers is bounded by
\begin{equation} \label{c3}
\mathop{\mathop{\sum\nolimits'}_{b \le z^{1/2}}}_{\text{gcd}(b,q) = 1} \mathop{\mathop{\sum_{x/b < n/b \le (x + y)/b}}_{p_{-}(n / b) > z^{1/2}}}_{(n/b) \equiv r \overline{b} \; (\bmod q)} 1 \ll \mathop{\sum\nolimits'}_{b \le z} \Bigl( \frac{y / b}{\phi(q) \log z} + z \Bigr) \ll \frac{y}{\phi(q) \log z} + z^{3/2} \ll_\alpha \frac{y}{\phi(q) \log y}
\end{equation}
by \eqref{asymptotic}, Lemma \ref{lem2} and the convergence of the sum of reciprocal of squarefull numbers (which follows from Lemma \ref{lem1} for instance). Here $\overline{b}$ denotes the multiplicative inverse of $b \;(\bmod{q})$, i.e. $b \overline{b} \equiv 1 \pmod{q}$.

\bigskip

Combining \eqref{c1}, \eqref{c2} and \eqref{c3}, we have Theorem \ref{thm2}.


\section{Proof of Theorem \ref{thm1.5}}

It is very similar to the proof of Theorem \ref{thm2}. So, we just highlight the necessary adjustments. We set $q = 1$ and $z = y^{1/3}$. The arguments for Case 1 and Case 2 are exactly the same as \eqref{c1} and \eqref{c2} and we get the bound
\[
\frac{y}{z^{1/4}} + z^{1/2} \ll y^{11/12}.
\]
It remains to deal with Case 3 where $b_n \le z^{1/2}$ and $z^{1/2} < p_{-}(d_n) \le y^{1/2}$ as the squarefull numbers are assumed to be $y^{1/2}$-smooth. Thus, with $p := p_{-}(d_n)$ and $d_n := p^2 d$, the number of squarefull numbers in this case is bounded by
\begin{align*}
\mathop{\sum\nolimits'}_{b \le z^{1/2}} \sum_{z^{1/2} < p \le y^{1/2}} \mathop{\sum_{x/b < p^2 d \le (x + y)/b}}_{p_{-}(d) > z^{1/2}} 1 =& \mathop{\sum\nolimits'}_{b \le z^{1/2}} \sum_{z^{1/2} < p \le y^{1/2}} \mathop{\sum_{\frac{x}{b p^2} < d \le \frac{x + y}{b p^2}}}_{p_{-}(d) > z^{1/2}} 1 \\
\ll& \mathop{\sum\nolimits'}_{b \le z^{1/2}} \sum_{z^{1/2} < p \le y^{1/2}} \Bigl( \frac{y/ (b p^2)}{\log z} + z \Bigr) \ll \frac{y}{z^{1/2} \log z} + \frac{z^{5/4} y^{1/2}}{\log y} \ll y^{11/12}
\end{align*}
by \eqref{asymptotic}, Lemma \ref{lem2} and the convergence of the sum of reciprocal of squarefull numbers. The above bounds together yield Theorem \ref{thm1.5}.


\section{Proof of Theorem \ref{thm2.5}}

Given integer $k \ge 3$ and small real number $\epsilon > 0$. We claim that, for some sufficiently small constant $c_\epsilon > 0$, the interval
\[
(x, x + c_\epsilon x^{1 - (2 + \epsilon)/k}]
\]
contains at most one $k$-full number for all sufficiently large $x \ge C$ (in terms of $\epsilon$ and $c_\epsilon$). We shall abbreviate $c_\epsilon x^{1 - (2 + \epsilon)/k}$ as $z$.

\bigskip

Suppose the contrary that the interval $(x, x + z]$ contains two $k$-full numbers $b < c$. Then $c = a + b$ for some integer $0 < a \le z$. Let $d = \text{gcd}(a, b)$. Then the integers $a/d$, $b/d$ and $c/d$ are pairwise relatively prime. Note that $\kappa(n) \le n^{1/k}$ for any $k$-full numbers. Apply the $abc$-conjecture to the equation $a/d + b/d = c/d$, we get
\begin{align*}
\frac{x}{d} < \frac{c}{d} &\ll_\epsilon \Bigl( \kappa \Bigl(\frac{a}{d} \frac{b}{d} \frac{c}{d} \Bigr) \Bigr)^{1 + \epsilon/2} = \Bigl( \kappa \Bigl(\frac{a}{d} \frac{b}{d} \Bigr) \kappa \Bigl(\frac{c}{d} \Bigr) \Bigr)^{1 + \epsilon/2} \\
&\le \Bigl( \kappa \Bigl(\frac{a b}{d} \Bigr) \kappa(c) \Bigr)^{1 + \epsilon/2} \le  \Bigl(\frac{a^{1/k} b^{1/k}}{d^{1/k}} c^{1/k} \Bigr)^{1 + \epsilon/2} \le \frac{z^{(1 + \epsilon/2)/k} (2x)^{2(1 + \epsilon/2) / k}}{d^{(1 + \epsilon/2)/k}}
\end{align*}
as $a b / d = \text{lcm}(a,b)$ is still $k$-full, $0 < a \le z$ and $b, c \le 2x$. Since $0 < d \le a \le z$, the above implies
\[
x^{1 - (2 +\epsilon)/k} \ll_\epsilon z
\]
which contradicts with $z = c_\epsilon x^{1 - (2 + \epsilon)/k}$ when $c_\epsilon > 0$ is small enough. Thus, the claim is true.

\bigskip

Clearly, Theorem \ref{thm2.5} is true for $1 \le y \le C$ by picking an appropriate implicit constant. Now, for $C \le y \le x$, the above claim implies that the interval
\[
(x, x + c_\epsilon y^{1 - (2 + \epsilon) / k}]
\]
contains at most one $k$-full number. By dividing the interval $(x, x+y]$ into subintervals of length $c_\epsilon y^{1 - (2 + \epsilon)/k}$, we obtain
\[
Q_k(x + y) - Q_k(x) \ll \frac{y}{c_\epsilon y^{1 - (2 + \epsilon) / k}} \cdot 1 \ll_\epsilon y^{\frac{2 + \epsilon}{k}}
\]
which gives Theorem \ref{thm2.5}.


\section{Proof of Theorem \ref{thm3}}

First, we suppose $y \le x^{0.2}$. We claim that there is no non-trivial three term arithmetic progression among the squarefull numbers in the interval $(x, x + y]$ under the abc-conjecture. Suppose the contrary. Then we have three squarefull numbers $x < a_1^2 b_1^3 < a_2^2 b_2^3 < a_3^2 b_3^2 \le x + y$ such that
\[
a_1^2 b_1^3 = a_2^2 b_2^3 - d \; \; \text{ and } \; \; a_3^2 b_3^3 = a_2^2 b_2^3 + d
\] 
for some positive integer $d$ with $2 d \le y$. Multiplying the above two equations, we get
\[
a_1^2 a_3^2 b_1^3 b_3^3 = a_2^4 b_2^6 - d^2 \; \; \text{ or } \; \; a_1^2 a_3^2 b_1^3 b_3^3 + d^2 = a_2^4 b_2^6.
\]
Say $D^2 = \text{gcd}(a_2^4 b_2^6, d^2)$ as the numbers are perfect squares. Then, the three integers
\[
\frac{a_1^2 a_3^2 b_1^3 b_3^3}{D^2}, \; \; \frac{d^2}{D^2}, \; \;  \frac{a_2^4 b_2^6}{D^2}
\]
are pairwise relatively prime and we have the equation
\[
\frac{a_1^2 a_3^2 b_1^3 b_3^3}{D^2} + \frac{d^2}{D^2} = \frac{a_2^4 b_2^6}{D^2}.
\]
Now, by the abc-conjecture, we have
\begin{align*}
\frac{x^2}{D^2} \le \frac{a_2^4 b_2^6}{D^2} &\ll_\epsilon \kappa \Bigl(\frac{a_1^2 a_3^2 b_1^3 b_3^3}{D^2} \frac{d^2}{D^2}  \frac{a_2^4 b_2^6}{D^2} \Bigr)^{1+\epsilon} \\
&\le \kappa \bigl( a_1^2 a_3^2 b_1^3 b_3^3 \bigr)^{1 + \epsilon} \kappa \Bigl(\frac{d^2}{D^2} \Bigr)^{1 + \epsilon} \kappa \bigl(a_2^4 b_2^6 \bigr)^{1 + \epsilon} \\
&\le (a_1 b_1 a_2 b_2 a_3 b_3)^{1 + \epsilon} \Bigl(\frac{d}{D}\Bigr)^{1 + \epsilon} \ll x^{3/2 + 3\epsilon/2} \frac{y^{1 + \epsilon}}{D^{1 + \epsilon}}.
\end{align*}
Since $1 \le D \le d \le y$, this implies $x^{1/2 - 3\epsilon/2} \ll_\epsilon D^{1 - \epsilon} y^{1+\epsilon} \ll y^2 \le x^{0.4}$ which is a contradiction for small enough $\epsilon$, say $\epsilon = 0.01$, and sufficiently large $x > C$ (in terms of the implicit constant).

\bigskip

Clearly, the theorem is true for $1 \le y \le C$ by picking an appropriate implicit constant. So, we may assume $y > C$. Since arithmetic progressions are invariant under translation, we may shift the interval $(x, x+y]$ to $(0, y]$. Therefore, by Theorem \ref{BSthm}, we have 
\[
Q_2(x + y) - Q_2(x) \ll \frac{y}{\log^{1+c} y}
\]
which gives the theorem.

\bigskip

Now, if $y > x^{0.2}$, one can simply divide the interval $(x, x + y]$ into subintervals of length $x^{0.2}$:
\[
(x, x + x^{0.2}] \cup (x + x^{0.2}, x + 2 x^{0.2}] \cup \cdots \cup \Bigl(x + \Bigl\lfloor \frac{y}{x^{0.2}}\Bigr\rfloor x^{0.2}, x + \Bigl( \Bigl\lfloor \frac{y}{x^{0.2}} \Bigr\rfloor + 1 \Bigr) x^{0.2} \Bigr]
\]
Then, over each interval $(x + i x^{0.2}, x + (i+1) x^{0.2}]$, we have the bound
\[
Q_2(x + (i+1) x^{0.2}) - Q_2(x + i x^{0.2}) \ll \frac{x^{0.2}}{\log^{1+c} x}.
\]
Summing over $\lfloor \frac{y}{x^{0.2}} \rfloor + 1$ of these intervals, we have
\[
Q_2(x + y) - Q_2(x) \ll \frac{y}{x^{0.2}} \cdot \frac{x^{0.2}}{\log^{1+c} x} \ll \frac{y}{\log^{1+c} y}
\]
which gives the theorem as well.

\bibliographystyle{amsplain}

Mathematics Department \\
Kennesaw State University \\
Marietta, GA 30060 \\
tchan4@kennesaw.edu

\end{document}